\documentclass[a4paper,12pt]{amsart}
\usepackage{amsmath,inputenc,euscript,amssymb,geometry}
\usepackage{graphicx}
\usepackage{amssymb}
\usepackage{latexsym}
\usepackage{amssymb,amsbsy,amsmath,amsfonts,amssymb,amscd}

\title[Dirac NLS]{On the Dirac delta as initial condition for nonlinear Schr\"odinger equations}

\thanks{*Second author supported by the grant MTM 2004-03029 of MEC 
(Spain) and FEDER. Both authors supported by the European Project HYKE (HPRN-CT-2002-00282).}
\thanks{
E-mail addresses : Valeria.Banica@univ-evry.fr, mtpvegol@lg.ehu.es}
\newtheorem{lemma}{Lemma}[section]
\newtheorem{theorem}[lemma]{Theorem}
\newtheorem{prop}[lemma]{Proposition}

\newtheorem{remark}[lemma]{Remark}

\newcommand{\mo}[1]{|#1|}

\newcommand{\RR}{\mathbb{R}}

\newcommand{\al}{\alpha}

\newcommand{\eps}{\epsilon}

\newcommand{\DP}[1]{{\partial_ #1}}

\newcommand{\D}{\displaystyle}

\def\tend{{\rightarrow}}
\newcommand{\R}{\mathbb{R}}

\begin{document}

\maketitle
\begin{center}{{\bf{V. Banica$^1$, L. Vega$^{2,*}$}}  \vspace{3mm}\\\tiny{$^1$D\'epartement de Math\'ematiques, Universit\'e d'Evry, France\\
$^2$Departamento de Matemáticas, Universidad del Pais Vasco, Spain }}
\end{center}

\bigskip

\bigskip
\noindent
------------------------------------------------------------------------------------------------------------------\\
{\small{{\bf{Abstract}}\\
In this article we will study the initial value problem for some Schr\"odinger equations with Dirac-like initial data and therefore with infinite $L^2$ mass, obtaining positive results for subcritical nonlinearities. In the critical case and in one dimension we prove that after some renormalization the corresponding solution has finite energy. This allows us to conclude a stability result in the defocusing setting. These problems are related to the existence of a singular dynamics for Schr\"odinger maps through the so called Hasimoto transformation.\\

\noindent
{\small{{\bf{R\'esum\'e}}\\
Dans cet article on \'etudie le probl\` eme de Cauchy pour des \'equations de Schr\"odinger avec donn\'ee initiale de type Dirac et donc avec masse $L^2$ infinie, obtenant des r\'esultats positifs pour les non lin\'earit\'es sous-critiques. Dans le cas critique et en dimension un, on montre qu'apr\`es une certaine renormalisation la solution correspondante est  d'\'energie finie. On en d\'eduit un r\'esultat de stabilit\'e dans le cas d\'efocalisant. Ces probl\`emes sont li\'es \`a l'existence d'une dynamique singuli\`ere des applications de type Schr\"odinger par la transformation de Hasimoto.
}}\\
\noindent
---------------------------------------------------------------------------------------------------------------------------\\
\bigskip
\section{Introduction}
 In this paper we  will study the IVP associated to the non-linear Schr\"odinger equation (NLS)
 \begin{equation}
\label{eq1}
\left\{\begin{array}{rcl}
iu_t+\Delta u\pm |u|^{\alpha}u&=&0\quad,\quad x\in\R^d\quad,\quad t>0\\
u(0,x)&=&a\delta_{x=0}+u_0(x)
\end{array}\right.
\end{equation}
with $\,0\leq\alpha$, and $\,u_0\,$ regular and possibly small with respect to $\,a\delta_{x=0}\,$. A related problem that we also consider is to take as $u(0,x)$ a perturbation of $ae^{i\frac{x^2}{4}}.$

Let us  recall first what is known about the broader question of  the minimal regularity needed to assume in the initial condition $u(0,x)$ so that \eqref{eq1} is well posed.
Within the Sobolev class
 $H^s$ with $s\geq 0$  the answer is positive and well understood, at least from the point of view of local wellposedness,  and is due to Ginibre and Velo \cite {GiVe}, Cazenave and Weissler  \cite {CaWe}. The proof follows a Picard iteration scheme based on  the so called Strichartz estimates. It was observed by Kenig, Ponce and Vega in \cite{KPV}, that if $s<0$  Picard's iteration can not work due to the lack of uniform continuity of the map datum-solution. Then Vargas and Vega propose in \cite{VV} a different class of spaces which are built using the Fourier transform of the initial condition. In the particular case of the one dimensional  cubic NLS they are able to consider a larger class than $L^2$. Finally Gr\"unrock in \cite{Gr} has extended that result, being able to prove local wellposednes  if  $u(0,x)$ satisfies that its Fourier transform is in some $L^p$ for $p<\infty$. Therefore he is just missing the delta function in \eqref{eq1} with $d=1$ and $\alpha=2$.

Let us recall next the known explicit solutions of  \eqref{eq1} if $\,u_0=0$. It was also noticed in \cite{KPV}, that if there is uniqueness of \eqref{eq1} with $\,u_0= 0\,$ then the corresponding solution should be invariant under the galilean transformations. That is to say for any $\,\nu\in\R^d$,
\begin{equation}
\label{eq2}
u_{\nu}(t,x):=e^{-it|\nu|^2+i\nu\cdot x}u(t, x-2\nu t)=u(t,x).
\end{equation}
This in turn implies that $\,u:=u_{a,\pm\alpha}$, with
\begin{equation}
\label{eq3}
u_{a,\pm\alpha}(t,x)=f_a(t,x)\,e^{\pm iA_{a,\alpha}(t)},
\end{equation}
where
\begin{equation}
\label{eq4}
f_a(t,x)=\D a\frac{e^{i\frac{|x|^2}{4t}}}{(it)^{d/2}},
\end{equation}
and
\begin{equation}
\label{eq5}
A_{a,\alpha}(t)=\left\{\begin{array}{rcl}
\D\frac{|a|^{\alpha}}{1-\alpha\frac d2}t^{1-\alpha\frac d2}\qquad \text{if}\quad \alpha\neq\frac{2}{d},\\
\\
|a|^{\frac{2}{d}}\, \log\,t\qquad \text{if}\quad \alpha=\frac{2}{d}.
\end{array}\right.
\end{equation}
Notice that
\begin{equation}
\label{eq6}
\lim_{t\downarrow 0}f_a(t,x)=a\delta_{x=0}\quad (\mathcal{S}^{\prime}),
\end{equation}
and as a conclusion the IVP \eqref{eq1} is ill-posed if $\,\alpha\ge \frac 2d$ --see Theorem 1.5 in \cite{KPV} for a precise statement.\\

A first natural question is if $\,u_{a,\pm\alpha}\,$ is a stable solution in the subcritical case $\,\alpha<\frac 2d$. We denote $\|f\|_2=\|f\|_{L^2}$.  We have the following result.

\begin{theorem}\label{Th1.1} Let $\,\alpha<\frac{2}{d}\,$. For $\,u_0\in L^2$, there exists a time $\,t_0=t_0(a,\|u_0\|_2)\,$ and a unique solution $\,u\,$ of \eqref{eq1} such that 
$$u-u_{a,\pm \alpha}\in L^p\left([0,t_0),L^q\right)\cap\mathcal{C}\left([0,t_0),L^2\right)$$
with $\,(p,q)\,$ any admissible pair,
$$\frac2p+\frac {d}{q}=\frac{d}{2}\,,\qquad 2\leq p\leq \infty\,,\quad\quad(d,p)\neq(2,2).
$$
Moreover if $\,0\leq\alpha\leq1\,$ then $\,t_0$ can be taken arbitrarily large.
\end{theorem}

The way to prove this theorem is to write
$$\eta(t,x)=e^{\mp iA_{\alpha,a}(t)}(u-u_{a,\pm\alpha}).$$ 
Using the fact that $f_a$ is a solution for the linear equation, we are lead to the equation
\begin{equation}
\label{eq7}
\left\{\begin{array}{rcl}
i\eta_t+\Delta\eta\pm\left(|\eta+f_a|^{\alpha}-|f_a|^{\alpha}\right)\left(\eta+f_a\right)&=&0\\
\eta(0,x)&=&u_0(x).
\end{array}\right.
\end{equation}
Then we see that if $\,\alpha<\frac{2}{d}\,$, the term $\,|f_a|^{\alpha}\,$ is locally integrable with respect to the time variable and therefore the usual Picard iteration scheme works. Moreover, if $\,0\leq\alpha\leq1\,$, we obtain an a priori control of the $\,L^2$ norm, so that we get a global result in this case.\\

Another natural question in the subcritical case is to consider a more regular perturbation of $\,u_{a,\pm\alpha}\,$. For this purpose, we introduce the conformal transformation
\begin{equation}
\label{conformal}
T(f)(t,x)=\D\frac{e^{i\frac{|x|^2}{4t}}}{(it)^{d/2}}f\left(\frac 1t,\frac xt\right).
\end{equation}
Let $\,w\,$ be defined by $u=Tw$. Then $\,u\,$ solves \eqref{eq1} for $\,0<t<t_0\,$ iff $\,w\,$ solves 
for $\,1/t_0<t<\infty$
$$-iw_t+\Delta w\pm\D\frac{1}{t^{2-\alpha\frac{d}{2}}}|w|^\alpha w=0.$$
That is iff $v(t,x)=w(t,x)e^{\mp iA_{a,\alpha}(1/t)}$ solves for $\,1/t_0<t<\infty$ the equation
\begin{equation}
\label{eqv}
-iv_t+\Delta v\pm\D\frac{1}{t^{2-\alpha\frac{d}{2}}}(|v|^\alpha-|a|^\alpha) v=0.
\end{equation}
Let us notice that  by the changes of variable we did, 
$$u(t,x)=T\left(e^{\pm iA_{a,\alpha}(1/\cdot)}\,v(\cdot,\cdot)\right)(t,x)=e^{\pm iA_{a,\alpha}(t)}T(v)(t,x),$$
the initial solution $u_{a,\pm\alpha}$ of \eqref{eq1} corresponds to 
the constant trivial solution $a$ of \eqref{eqv}. 
Finally, the equation of the perturbation $\epsilon(t,x)=v(t,x)-a$, with initial data at time $1/t_0$, writes
\begin{equation}
\label{eqeps}
\left\{\begin{array}{rcl}
-i\epsilon_t+\Delta \epsilon\pm\D\frac{1}{t^{2-\alpha\frac{d}{2}}}(|\epsilon+a|^\alpha-|a|^\alpha) (\epsilon+a)&=&0,\\
\epsilon(1/t_0,x)=\epsilon_0(x).
\end{array}\right.
\end{equation}
We shall study this equation for large times, in appropriate Sobolev spaces, and under suitable conditions on $\alpha$. The subcritical condition $\alpha<\frac{2}{d}$ will be crucial in the proofs, since it gets the integrability at infinity of the time-coefficient in the nonlinearity. The asymptotic behavior of $\epsilon$ will give us informations on equation  \eqref{eq1}, and more precisely on the small time behavior of the perturbations around the solution $u_{a,\pm\alpha}$. 
The fact that we have been able to prove Theorem 1.1 directly on the initial equation \eqref{eq1} is related to the fact that the mixed spaces we are working with are invariant by the conformal transformation. \par
We have the following result.

\begin{theorem}\label{subcrit}  
Let $\,\alpha<\frac 2d\,$,  and let $s>\frac{d}{2}$. For $\epsilon_0\in H^s$, with norm small with respect to $|a|$, there exists a time $t_0=t_0(a)$ such that equation \eqref{eqeps} has a unique solution $\epsilon$ in a small ball of $L^\infty((1/t_0,\infty),H^s)$. Moreover, the wave operator exists and the equation enjoys the property of asymptotic completeness in $H^s$.\\
As a consequence, for $u_0$ small in $\Sigma^s=\{(1+|x|^s)f\in L^2\}$ with respect to $|a|$ , there exists a time $t_0=t_0(a)$ such that equation \eqref{eq1} admits a solution that writes, for all $0<t<t_0$,
\begin{equation}
\label{Hs}
u(t,x)=u_{a,\pm\alpha}+\D a\frac{e^{i\frac{|x|^2}{4t}}}{(it)^{d/2}}e^{\pm i A_{a,\alpha}(t)}\epsilon\left(\frac{x}{t},\frac{1}{t}\right),
\end{equation}
for a unique $\epsilon$ small in $L^\infty((1/t_0,\infty),H^s)$.
\end{theorem}

The proof of the results on $\epsilon$ is again standard and relies on the fact that in the case $s>\frac{d}{2}$, $H^s$ is an algebra included in $L^\infty$. As usual there is the difficulty of the lack of regularity of the nonlinear function appearing in \eqref{eqeps}, but this is overcome assuming smallness of $\epsilon$ with respect to $a$. The passage back to the initial equation \eqref{eq1} is done by using the scattering results combined with the asymptotic behavior of the linear Schr\"odinger evolution.\\

A case of rougher perturbations of $a\delta_0$ was recently studied by Kita. In \cite{Ki} he described the structure of the solutions yielded by initial data exactly a sum of two or three Dirac masses.

Finally, let us say a few words about which is the situation in the parabolic setting of the IVP analogous to  \eqref{eq1} , that is to say
\begin{equation}\label{eq111}
\left\{\begin{array}{rcl}
u_t-\Delta u\pm |u|^{\alpha}u&=&0\qquad x\in\R^d\quad,\quad t>0\\
u(0,x)&=&a\delta_{x=0}.\end{array}\right.
\end{equation}
This equation has been  intensively studied. It particular it has been shown by Weissler \cite{We} that in the focusing case (sign $-$), for $\al\leq \frac 2d\,$, there is no uniqueness. For the defocusing case, if only positive solutions are considered, for $\al<\frac 2d\,$ there is a unique solution, and for $\al\geq \frac 2d\,$ there is no such solution. This result was proved by Br\'ezis and Friedman in \cite{BrFr}. So, for the heat equation, an important difference is made between the focusing and the defocusing case, even for $\al<\frac 2d\,$. This is not the picture for the Schr\"odinger equation as we see from the statements of Theorem \ref{Th1.1}  and Theorem \ref{subcrit}.\\

\bigskip

Let us consider next the critical case $\,\alpha=\frac 2d\,$ and just in dimension $\,d=1$. As it is well known this equation is completely integrable. It is also closely related to the Schr\"odinger map equation
\begin{equation}
\label{eq8}
\left\{
\begin{array}{rcl}
\gamma_t&=&\gamma\land_{\pm}\gamma_{xx}\\
\gamma(0,x)&=&\gamma_0
\end{array}\right.
\end{equation}
with $\,b\land_{\pm}c:=\mathcal{A}_{\pm}(b\land c)$, and
$$\mathcal {A}_{\pm}=\left(\begin{array}{ccr}
1&0&0\\
0&1&0\\
0&0&\pm1
\end{array}\right).$$
Then it is straightforward that
\begin{equation}
\label{eq10}
\DP{t}\,\langle {\mathcal {A}}_{\pm}\gamma,\gamma\rangle=0
\end{equation}
if $\,\gamma\,$ solves \eqref{eq8}. Therefore for $\,\mathcal {A}_+\,$ we get the Sch\"odinger map onto the unit sphere $\mathbb S^2$, and for $\,\mathcal {A}_-\,$ we have it onto 2d hyperbolic space:
$$\mathbb {H}^2=\{\,(a_1,a_2,a_3)\in\R^3\quad \text{such that} \quad a_1^2+a_2^2-a_3^2=-1,\qquad a_3>0\,\}.
$$

Equation \eqref{eq8} can be also obtained from the flows of curves in $\,\R^3\,$  given by
\begin{equation}
\label{eq11}
\chi_t=\chi_x\land_{\pm}\chi_{xx}.
\end{equation}

In the euclidean case this equation is also called the Local Induction Approximation and was obtained by Da Rios \cite{DaR} as a crude model which describes the dynamics of a vortex filament equation within Euler equations. In \cite{SJL}, Guti\'errez, Rivas and Vega obtain solutions of the IVP
\begin{equation}
\label{eq12}\left\{
\begin{array}{rcl}
\chi_t&=&\chi_x\land_{+}\chi_{xx},\\\\
\chi(0,x)&=&{\left\{\begin{array}{ll}
A_1^{+}x&x\ge 0,\\
A_2^{+}x&x\le 0,
\end{array}\right.}
\end{array}\right.
\end{equation}
for any unit vectors $\,A_1^{+}\,,\,A_2^{+}\,$ such that
$$A_1^{+}+A_2^{+}\neq0.$$

More recently  de la Hoz \cite{Pa} has proved a similar result in the non--elliptic setting. Namely he obtains solutions of
\begin{equation}
\label{eq13}\left\{
\begin{array}{rcl}
\chi_t&=&\chi_x\land_{-}\chi_{xx},\\\\
\chi(0,x)&=&{\left\{\begin{array}{ll}
A_1^{-}x&x\ge 0,\\
A_2^{-}x&x\le 0,
\end{array}\right.}
\end{array}\right.
\end{equation}
for any pair of vectors $A^-_1$ and $A^-_2$ in $\mathbb H^2$.
 
In both cases the corresponding solution curves are described geometrically by the curvature and the (generalized) torsion given by
\begin{equation}
\label{eq15}
c(t,x)=\D\frac{c_0}{\sqrt t}\quad,\quad \tau(t,x)=\D\frac x{2t},
\end{equation}
with $\,c_0\,$ a free parameter uniquely determined by $\,\left(A_1^\pm,A_2^\pm\right)$. Therefore at time $\,t=1\,$ the curves are real analytic while at $\,t=0\,$ a corner is developed if $c_0\neq0$. Recall that the flows given in \eqref{eq11} are both reversible in time so that the above solutions are examples of the formation of a singularity in finite time. From \eqref{eq15} we conclude that for all time
$$\int_{-\infty}^\infty\,c^2(t,x)\,dx=\infty.
$$
However we shall see below that after some renormalization these curves  have finite energy.

The connection of \eqref{eq11} with the 1d cubic NLS was proved by Hasimoto who used the transformation
\begin{equation}
\label{eq16}
\Psi(t,x)=c(t,x)\exp{\left\{i\int\limits_0^x \tau(t,x')dx'\right\}}.
\end{equation}

Although he worked just in the euclidean case (i.e. $\,\land_+\,$ in \eqref{eq11}), a similar argument can be given for the other case, see for example \cite{Ko} and \cite{Di}. The final conclusion is that if $\,\chi\,$ solves \eqref{eq11} and $\,\Psi\,$ is defined as in \eqref{eq16} then
\begin{equation}
\label{eq17}
i\Psi_t+\Psi_{xx}\pm\D\frac 12\left(|\Psi|^2+a(t)\right)\Psi=0,
\end{equation}
for some real function $a(t).$
So in the particular case given in \eqref{eq15} we have
$$\Psi(t,x)=c_0\D\frac{e^{i \frac{x^2}{4t}}}{\sqrt t}=\sqrt{i}f_{c_0}(t,x).$$
Then chosen $\,a(t)=-\D\frac{|c_0|^2}{t}\,$ we get a solution of \eqref{eq17} with
$$\Psi(0,x)=\sqrt i\,c_0\delta.$$
Notice that the coefficient $\,\D\frac 12\,$ which appears in \eqref{eq17} is harmless and can be easily absorbed by the change of variable
$$(t',x')=\left(2t, \sqrt 2 x\right).$$

Hence we are interested in solving on $\mathbb{R}$, around the particular solution $f_a$, the equation
\begin{equation}
\label{eq18}
\left\{\begin{array}{rcl}
iu_t+u_{xx}\pm\left(|u|^2-\D\frac{|a|^2}{t}\right)u&=&0,\\
u(0,x)&=&a\delta_{x=0}+u_0(x).
\end{array}\right.
\end{equation}
Unfortunately we are not able to deal with \eqref{eq18} directly, so that we propose for any $t_0>0$ the related problem 
\begin{equation}
\label{eq19}
\left\{\begin{array}{rcl}
iu_t+u_{xx}\pm\left(|u|^2-\D\frac{|a|^2}{t}\right)u&=&0,\\
u(t_0,x)&=&a\D\frac{e^{i\frac{x^2}{4t_0}}}{\sqrt{it_0}}+u_1(x),
\end{array}\right.
\end{equation}
and we look for a backward solution. That is to say for
\begin{equation}
\label{eq19prima}
0<t<t_0.
\end{equation}
It is natural to consider as before the conformal transformation
\begin{equation}
\label{eq192prima}
u(t,x)=Tv(t,x)=\D\frac{e^{i\frac{x^2}{4t}}}{\sqrt{it}}v\left(\frac 1t,\frac xt\right).
\end{equation}
Then $\,u\,$ solves \eqref{eq19} for $\,0<t<t_0\,$ iff $\,v\,$ solves for $ 1/t_0<t<\infty$ 
\begin{equation}
\label{eq20}
\left\{\begin{array}{rcl}
-iv_t+v_{xx}\pm\D\frac 1t\left(|v|^2-|a|^2\right)v&=&0,\\
v\left(1/t_0,x\right)&=&a+\epsilon_0(x),
\end{array}\right.
\end{equation}
where $\,v_0\,$ is defined by $\,u_1(x)=(T\epsilon_0)(t_0,x)$. It is easy to solve \eqref{eq20} locally in time for both situations focussing and defocussing. Moreover  there is a natural energy. In fact if we define
\begin{equation}
\label{eq21}
E(t)=\D\frac 12\int|v_x(t,x)|^2dx\mp\frac 1{4t}\int\left(|v(t,x)|^2-|a|^2\right)^2dx,
\end{equation}
then if $\,v\,$ is a solution of \eqref{eq20}, we have
\begin{equation}
\label{eq22}
\DP{t}E(t)\mp\frac 1{4t^2}\int\left(|v(t,x)|^2-|a|^2\right)^2dx=0.
\end{equation}

As a consequence, and in the defocussing situation, we will be able to prove the following theorem which is the main result of this paper.

\begin{theorem}\label{Th1.3} For all $\,t_0>0\,$ and for all $\,\epsilon_0\in H^1$, there exists a unique solution of the IVP
\begin{equation}
\label{eq23}
\left\{\begin{array}{rcl}
-iv_t+v_{xx}-\D\frac 1t\left(|v|^2-|a|^2\right)v&=&0,\quad 1/t_0<t<\infty,\\
v\left(1/t_0,x\right)&=&a+\epsilon_0(x),
\end{array}\right.
\end{equation}
with 
$$v-a\in\mathcal{C}\left((1/t_0,\infty),H^1\right).$$
Moreover
\begin{equation}
\label{eq24}
\int|v_x(t,x)|^2dx\le 2E(1/t_0),
\end{equation}
\begin{equation}
\label{eq25}
\int_{\frac{1}{t_{0}}}^{\infty}\int\left(|v(t,x)|^2-|a|^2\right)^2dx\D\frac{dt}{t^2}<4E(1/t_0),
\end{equation}
and in particular
\begin{equation}
\label{eq26}
\liminf_{t\to\infty}\D\frac 1t\int\left(|v(t,x)|^2-|a|^2\right)^2dx=0.
\end{equation}
\end{theorem}

\begin{remark} Using \eqref{eq26} and \eqref{eq192prima}
 we obtain that the solution $u$ of \eqref{eq19} in the defocussing setting for
$$0<t<t_0$$
satisfies
\begin{equation}
\label{eq27}
\liminf_{t\to 0}\left\|t|u(t)|^2-|a|^2\right\|_2=0,
\end{equation}
which can be understood as a weak stability result of the singular solution $\,a\D\frac{e^{i\frac {x^2}{4t}}}{\sqrt t}\,$ of \eqref{eq19}. 
\end{remark}

\begin{remark} 
The limit of $u(t,x)$ when $t$ goes to zero is not settled                                  
with our approach. In a forthcoming paper                                       
we shall look at the problem as a long range                                    
scattering one assuming some appropiate                                         
smallness condition. \\

\end{remark}

\begin{remark} It is interesting to write which is the energy for \eqref{eq16}, the solution  of \eqref{eq17}, in terms of the geometric quantities $c$ and $\tau$. It is given by
\begin{eqnarray}
\nonumber\widetilde E(t) & = & {t^2\over
4\sqrt2}\int_{-\infty}^{+\infty}\bigg(c_x^2(t, x) + c^2(t,
x)\left({x\over 2t} - \tau(t, x)\right)^2\bigg)
dx \\
& + & {1\over 16\sqrt 2}\int_{-\infty}^{+\infty}[t c^2(t, x) -
c_0^2]^2dx.
\end{eqnarray}

\noindent Then
\begin{equation}
{d\over dt}\widetilde E(t) - {1\over 16\sqrt 2\ t}
\int_{-\infty}^{+\infty} [t c^2(t, x) - c_0^2]^2dx = 0,
\end{equation}
and
$$
\liminf_{t\to 0}\left\|t|c|^2-|c_0|^2\right\|_2=0.
$$
Recall that in this case we will solve the equation backwards in time.

\end{remark}

Finally, let us make a remark on the Gross-Pitaevskii defocusing equation
\begin{equation}\label{GP}
\left\{\begin{array}{rcl}
i\psi_t+\Delta \psi-\left(|\psi|^2-1\right)\psi&=&0,\\
\psi(0,x)&=&\psi_0(x).
\end{array}\right.
\end{equation}
This equation was globally solved in $1+H^1(\mathbb{R}^d)$, with an exponential growth in time control of the mass $\|\psi(t) -1\|_2$, for $d\in\{ 2,3\}$ by Bethuel and Saut, and for $d=1$ by Gallo (\cite{BS},\cite{G2}).\par
In a class of larger spaces, the Zhidkov spaces $X^k(\mathbb{R}^d)$, it has been solved in $X^1(\RR)$ and in $X^2(\RR^2)$ by Zhidkov, by Gallo and by Goubet (\cite{Z1}, \cite{Z2}, \cite{G1},\cite{Go}). 
Also, considered in the natural energy space $\{f\in H^1_{loc} \,\,,\,\,\nabla f\in L^2\,\,,\,\,|f|^2-1\in L^2\}$, it has been solved for $d\in\{ 2,3\}$ and for $d=4$ with smallness assumption, by G\'erard in \cite{P}. Recently, Gustafson, Nakanishi and Tsai described in \cite{GNT} the scattering in modified  $1+H^{\frac{d}{2}-1}(\mathbb{R}^d)$ spaces, for $d\geq 4$ and small data. The proof relies on the linearized equation and the conservation of the energy is not used.\par
By exploiting more the mass and energy laws, as done for proving Theorem \ref{Th1.3}, we can get a slightly modified proof of the very short and simple one in \cite{BS}, allowing to have the following result.
\begin{prop}\label{propGP}
The solution of (\ref{GP}) is globally well-posed in $1+H^1(\mathbb{R}^d)$, for all dimensions $d$ such that local existence occurs, that is surely for $d\in\{1,2,3\}$, with the control
$$\|\psi(t) -1\|_2\leq c\, t,$$
where the constant $c$ depending on the initial data.
\end{prop}

The article is structured as follows. The first two sections contain the proofs for the results in the subcritical case. In section \S 4 we prove Theorem \ref{Th1.3}. Section \S5 concerns Proposition \ref{propGP}, and in the last section we give a technical lemma. \\

\noindent
{\bf{Acknowledgements.}} We thank the referee for the careful reading of the paper and his multiple suggestions.



\section{The sub-critical power. The Strichartz case}

We shall prove the following Lemma, that implies the local existence result in Theorem \ref{Th1.1}, as indicated in the introduction.
\begin{lemma} \label{glestr}
There exists $t_0=t_0(a,\|u_0\|_2)$ and a unique solution $\eta$ of the equation (\ref{eq7})
$$\left\{\begin{array}{rcl}
i\eta_t+\Delta \eta\pm (\mo{\eta+f_a}^\al-|f_a|^\al)(\eta+f_a)&=&0,\\
\eta(0,x)&=&u_0(x),
\end{array}\right.$$
such that
$$\eta\in L^\infty([0,t_0),L^2)\cap L^p([0,t_0),L^q),$$
with $(p,q)$ any admissible couple
$$\frac2p+\frac {d}{q}=\frac{d}{2}\,,\quad\quad2\leq p\leq \infty\,,\quad\quad(d,p)\neq(2,2).
$$

\end{lemma}

\begin{proof} 
Let us denote $X$ the intersection of the mixed spaces. In order to do a fixed point argument in a closed ball of $X$, we have to estimate the norm of the operator
$$\Phi(\eta)(t,x)=e^{it\Delta}u_0(x)\pm i\int_0^t e^{i(t-\tau)\Delta}\left(|\eta(\tau,x)+f_a(\tau,x)|^\al-|f_a(\tau,x)|^\al\right)(\eta(\tau,x)+f_a(\tau,x))d\tau.$$
The Schr\"odinger operator is unitary on $L^2$, so one gets 
\begin{equation}\label{str}
\|\Phi(\eta)(t)\|_2\leq \|u_0\|_2+\int_0^t \left\|   \left(|\eta(\tau)+f_a(\tau)|^\al-|f_a(\tau)|^\al\right)(\eta(\tau)+f_a(\tau)) \right\|_2d\tau.\end{equation}
The nonlinearity we are working with,
$$F(z)=(|z+f_a(\tau,x)|^\alpha-|f_a(\tau,x)|^\alpha)(z+f_a(\tau,x)),$$ 
verifies $F(0)=0$ and
$$|F'(z)|=\max\{|\partial_z F(z)|,|\partial_{\overline{z}}F(z)|\}\leq c\,(|f_a(\tau,x)|^\al+|z+f_a(\tau,x)|^\alpha).$$
Under the assumption $\alpha\geq 0$, we get
$$|F'(z)|\leq c\,(|f_a(\tau)|^\alpha+|z|^\alpha).$$
This is the classical growth hypothesis on nonlinearities for proving the local wellposedness. As done for the local $L^2$ Cauchy problem, one can split the nonlinearity in two parts, $F=F_1+F_2$, with $|F_1'(z)|\leq c\,|f_a(\tau)|^\alpha$ and $|F_2'(z)|\leq c\,|z|^\alpha$, use 
$$F(z_1)-F(z_2)=\int_0^1F'(tz_1+(1-t)z_2)\,dt \,(z_1-z_2),$$
for estimating $F(\eta(\tau,x))-F(0)=F(\eta(\tau,x))$, and obtain 
$$\|\Phi(\eta)(t)\|_2\leq \|u_0\|_2+C\int_0^t |f_a(\tau)|^\al\|\eta(\tau)\|_2+\|\eta^{\al+1}(\tau)\|_2\,d\tau.$$
Therefore, for all $t>0$,
$$\|\Phi(\eta)(t)\|_2\leq \|u_0\|_2+c(a)\|\eta\|_X\,\int_0^t\frac{d\tau}{\tau^{\alpha\frac{d}{2}}}+\int_0^t \|\eta(\tau)\|_{L^{2(\al+1)}}^{\al+1}d\tau.$$
Now, we compute the second integral, and in the last term we perform a H\"older inequality
$$\|\Phi(\eta)(t)\|_2\leq \|u_0\|_2+c(a)\|\eta\|_X \,\,t^{1-\alpha\frac{d}{2}}+c\|\eta\|_{L^p([0,t],L^q)}^{\al+1}\,t^\frac{4-d\alpha}{4},$$
where
$$p=\frac{4(\al+1)}{d\al}\:\:\:\:q=2(\al+1),$$
form an admissible couple. 
We get 
$$\underset{0<t\leq t_0}{\sup}\|\Phi(\eta)(t)\|_2\leq \|u_0\|_2+c(a)t_0^{1-\alpha\frac{d}{2}}\|\eta\|_X+ct_0^{1-\alpha\frac{d}{4}}\|\eta\|_X^{\al+1}.$$
Then, for $\|u_0\|_2$ finite, $\Phi(\eta)$ satisfies the first condition to be in $X$. \\

Let us treat now the $L^pL^q$ norm. By using the inhomogeneous global Strichartz inequalities, we obtain that for $t>0$,
$$\|\Phi(\eta)(t)\|_{L^pL^q}\leq \|u_0\|_2+\int_0^{t_0} \left\|\left(|\eta(\tau)+f_a(\tau)|^\al-|f_a(\tau)|^\al\right)(\eta(\tau)+f_a(\tau))\right\|_2d\tau.$$
The right hand side term can be treated exactly like the one in (\ref{str}). In conclusion, if $u_0$ is in $L^2$, $\Phi(\eta)$ stays in $X$. Arguing as before, and assuming $t_0$ small with respect to $\|u_0\|_2$ and to $|a|$, we get that the operator $\Phi$ is a contraction on a closed ball of $X$. 
Therefore the lemma follows from the fixed point theorem.
\end{proof}

Let us prove now the global existence result of Theorem \ref{Th1.1}, in the case $\al\in[0,1]$. 
By multiplying equation (\ref{eq7}) by $\overline{\eta}$, and by taking then the imaginary part, we have
\begin{equation}\label{massteta}
\DP{t}\|\eta(t) \|_2^2=\mp \,2\,\int \left(|\eta(t,x)+f_a(t,x)|^{\al}-|f_a(t,x)|^{\al}\right)\Im f_a(t,x)\overline{\eta}(t,x)dx.
\end{equation}

\begin{lemma}\label{masscontroleta}
If $\al\in[0,1]$ there exists $c(a)>0$ such that for all $t\geq t'$ we have the a priori estimate 
$$\|\eta(t)\|_2^2\leq \|\eta(t')\|_2^2 \,e^{c(a)t^{1-\alpha\frac{d}{2}}}.$$
\end{lemma}

\begin{proof}
Lemma \ref{tech} allows us to upper-bound the right hand side of (\ref{massteta}) and get
$$\DP{t}\|\eta (t)\|_2^2\leq\frac{c(a)}{t^{\alpha\frac{d}{2}}}\|\eta (t)\|_2^2,$$
and the lemma follows.
\end{proof}

Therefore, if $\al\in[0,1]$, by using Lemma \ref{masscontroleta} we can extend the solution $\eta$ for arbitrary large $t_0$, as done in the critical case in \S 4.2. 


\section{The sub-critical power. The $H^s$ case}
\subsection{Existence of solutions}

\begin{lemma}\label{glexistsubcrit}
Let $s>\frac{d}{2}$, and let $\eps_0\in H^s$ with
$$\|\eps_0\|_{H^s}\leq \frac{|a|}{8}.$$
There exists a time $t_0=t_0(a)$ such that equation (\ref{eqeps}),
$$\left\{\begin{array}{rcl}
-i\epsilon_t+\Delta \epsilon\pm\D\frac{1}{t^{2-\alpha\frac{d}{2}}}(|\epsilon+a|^\alpha-|a|^\alpha) (\epsilon+a)&=&0,\\
\epsilon(1/t_0,x)=\epsilon_0(x),
\end{array}\right.$$
has a unique solution $\epsilon$ in
$$Y=\left\{f\in L^\infty((1/t_0,\infty),H^s) , \underset{t\geq 1/t_0}{\sup}\,\|f(t)\|_{H^s}\leq\frac{|a|}{4}\right\}.$$
\end{lemma}

\begin{proof}
 
In order to do a fixed point argument in this space, we have to estimate the norm in $Y$ of the operator
$$\Phi(\eps)(x,t)=e^{-it\Delta}\eps_0(x)\pm i\int_{1/t_0}^t e^{-i(t-\tau)\Delta}\left(|\eps(x,\tau)+a|^\al-|a|^\al\right)(\eps(x,\tau)+a)\frac{d\tau}{\tau^{2-\al\frac{d}{2}}}.$$
By using the fact that the Schr\"odinger operator is unitary on $H^s$, 
$$\|\Phi(\eps)(t)\|_{H^s}\leq \|\eps_0\|_{H^s}+
\int_{1/t_0}^t \|\left(|\eps(\tau)+a|^\al-|a|^\al\right)(\eps(\tau)+a)\|_{H^s}\frac{d\tau}{\tau^{2-\al\frac{d}{2}}} .$$
The nonlinearity
$$\tilde{F}(z)=(|z+a|^\alpha-|a|^\alpha)(z+a)$$
is a $\mathcal{C}^\infty$ function on $|z|<\frac{|a|}{4}$ with $\tilde{F}(0)=0$. Here we see the difference with the case of classical power-nonlinearity $|z|^\alpha z$, whose lack of regularity imposes the conditions $s<\alpha$ or $\alpha$ even. 
Since $\epsilon(\tau)$ is in $Y$ and $H^s$ is embedded in $L^\infty$, it follows that $|\eps(\tau)|\leq\|\epsilon(\tau)\|_{H^s}<\frac{|a|}{4}$. Therefore (\cite{CM})
$$\|\tilde{F}(\epsilon(\tau))\|_{H^s}\leq c(a),$$
and
$$\|\Phi(\eps)(t)\|_{H^s}\leq \|\eps_0\|_{H^s}+c(a)\int_{1/t_0}^t \frac{d\tau}{\tau^{2-\al\frac{d}{2}}}.$$
By computing the integral,
$$\|\Phi(\eps)(t)\|_{H^s}\leq \|\eps_0\|_{H^s}+c(a)\frac{|t^{-1+\al\frac{d}{2}}-t_0^{1-\al\frac{d}{2}}|}{|-1+\al\frac{d}{2}|}.$$
We are in the case $\al<\frac{2}{d}$, so 
$$\|\Phi(\eps)\|_Y\leq \|\eps _0\|_{H^s}+c(a)t_0^{1-\al\frac{d}{2}}.$$
The smallness of $\eps_0$ in $H^s$ in $Y$ yield
$$\|\Phi(\eps)\|_Y\leq \frac{|a|}{8}+c(a)t_0^{1-\al\frac{d}{2}}.$$
Then, for $t_0=t_0(a)$ small enough, $\Phi(\eps)$ remains in $Y$. Moreover, by arguing similarly, we obtain that the operator acts as a contraction on $Y$. Therefore the fixed point theorem ends the proof.
\end{proof}

\subsection{Scattering properties}

\begin{lemma}\label{waveop} Let $s$ an integer such that $s>\frac{d}{2}$.\\
i) For all $\eps_+\in H^s$ with
$$\|\eps_+\|_{H^s}\leq \frac{|a|}{8},$$
there exists a solution $\eps$ of (\ref{eqeps}) in $Y$ such that
$$\underset{t\tend +\infty}{\lim}\|\eps(t)-e^{-it\Delta}\eps_+\|_{H^s}=0.$$
ii) Let $\eps_0\in H^s$ with
$$\|\eps_0\|_{H^s}\leq \frac{|a|}{8},$$
and let $\eps$ the solution given by Lemma \ref{glexistsubcrit}. Then there exists a $\eps_+$ small in $H^s$ such that 
$$\underset{t\tend +\infty}{\lim}\|\eps(t)-e^{-it\Delta}\eps_+\|_{H^s}=0.$$
\end{lemma}

The wave operator of i) is obtained by doing the same calculus as in the previous subsection, for a fixed point in $Y$ with the operator 
$$\Phi(\eps)(x,t)=e^{-it\Delta}\eps_+(x)\pm i\int_{t}^{+\infty} e^{-i(t-\tau)\Delta}\left(|\eps(x,\tau)+a|^\al-|a|^\al\right)(\eps(x,\tau)+a)\frac{d\tau}{\tau^{2-\al\frac{d}{2}}}.$$
Since the Schr\"odinger operator is unitary on $H^s$, proving ii) is equivalent to proving that 
$$e^{it\Delta}\eps(t)$$
has a limit in $H^s$ as $t$ goes to infinity, and also that
$$\lim_{t_1,t_2\tend +\infty}\| e^{it_1\Delta}\eps(t_1)-e^{it_2\Delta}\eps(t_2)\|_{H^s}=0.$$
By using the Duhamel formulation, it is enough to show 
$$\lim_{t_1,t_2\tend +\infty}\int_{t_1}^{t_2}\|\left(|\eps(\tau)+a|^\al-|a|^\al\right)(\eps(\tau)+a)\|_{H^s}\frac{d\tau}{\tau^{2-\al\frac{d}{2}}}=0.$$
The last assertion follows from calculus similar to the one in the proof of Lemma \ref{glexistsubcrit}.

\subsection{Proof of Theorem \ref{subcrit}}
For $u_0$ small in $\Sigma=\{(1+|x|^s)f\in L^2\}$, let us define $\eps_+$ by the Fourier relation, $\hat{\eps_+}(x/2)=u_0(x)$. Then $\eps_+$ is small in $H^s$ and Lemma \ref{waveop} ensures us that there exists a time $t_0=t_0(a)$ and a unique solution $\eps$ of (\ref{eqeps}) in $L^\infty((1/t_0,\infty),H^s)$ behaving at infinity as the free evolution of $\eps_+$.\\
Therefore we have the existence of a solution of the equation
$$i\DP{t}u+\Delta u\pm|u|^\al u=0,$$
for all $0<t<t_0$, given by
\begin{equation}\label{Teps}
u(x,t)=u_{a,\pm\alpha}+\frac{e^{i\frac{|x|^2}{4t}}}{(it)^\frac{d}{2}}e^{\pm iA_{a,\alpha}(t)}\eps\left(\frac{1}{t},\frac{x}{t}\right).
\end{equation}
Let us see now what happens with $u$ at time $0$. 
On one hand,
$$\left\|\frac{e^{i\frac{|x|^2}{4t}}}{(it)^\frac{d}{2}}e^{\pm iA_{a,\alpha}(t)}\left(\eps\left(\frac{1}{t},\frac{x}{t}\right)-e^{-\frac{i}{t}\Delta}\eps_+\left(\frac{x}{t}\right)\right)\right\|_2\leq \frac{1}{t^\frac{d}{2}} \left\|\eps\left(\frac{1}{t},\frac{x}{t}\right)-e^{-\frac{i}{t}\Delta}\eps_+\left(\frac{x}{t}\right)\right\|_2=$$
$$= \left\|\eps\left(\frac{1}{t}\right)-e^{-\frac{i}{t}\Delta}\eps_+\right\|_2.$$
Then the decay of Lemma \ref{waveop} allows us to say that 
\begin{equation}\label{uno}
\lim_{t\downarrow 0}\left\|\frac{e^{i\frac{|x|^2}{4t}}}{(it)^\frac{d}{2}}e^{\pm iA_{a,\alpha}(t)}\left(\eps\left(\frac{1}{t},\frac{x}{t}\right)-e^{-\frac{i}{t}\Delta}\eps_+\left(\frac{x}{t}\right)\right)\right\|_2=0.
\end{equation}
On the other hand, using the free Schr\"odinger evolution,
$$\frac{e^{i\frac{|x|^2}{4t}}}{(it)^\frac{d}{2}}e^{\pm iA_{a,\alpha}(t)}e^{-\frac{i}{t}\Delta}\eps_+\left(\frac{x}{t}\right)=
\frac{e^{i\frac{|x|^2}{4t}}}{(it)^{\frac{d}{2}}}e^{\pm iA_{a,\alpha}(t)}c\frac{e^{-i\frac{|x|^2}{4t}}}{(-\frac{i}{t})^{\frac{d}{2}}}\int e^{-iy^2\frac{t}{4}}e^{i\frac{xy}{2}}\eps_+(y) dy.$$
So, taking in account $\lim_{t\downarrow 0}A_{a,\alpha}(t)=0$, we have  
\begin{equation}\label{due}
\lim_{t\downarrow 0}\left\|\frac{e^{i\frac{|x|^2}{4t}}}{(it)^\frac{d}{2}}e^{\pm iA_{a,\alpha}(t)}e^{-\frac{i}{t}\Delta}\eps_+\left(\frac{x}{t}\right)-\hat{\eps}_+\left(\frac{x}{2}\right)\right\|_2=0.
\end{equation}
Therefore, in view of (\ref{Teps}), (\ref{uno}) and (\ref{due}), we conclude that $u$ verifies the initial condition
$$u(0,x)=a\delta_{x=0}+u_0(x),$$
so that $u$ is a solution of equation (\ref{eq1}).

\section{The critical power, defocusing case}\label{defocusing}
In order to treat the defocusing equation (\ref{eq23}), for $1/t_0<t<\infty$,
$$\left\{\begin{array}{rcl}
-iv_t+v_{xx}-\D\frac 1t\left(|v|^2-|a|^2\right)v&=&0,\\
v\left(1/t_0,x\right)&=&a+\epsilon_0(x),
\end{array}\right.$$
we shall consider instead the equation on $\eps=v-a$,  
\begin{equation}\label{ee}\left\{\begin{array}{rcl}
-i\eps_t+\eps_{xx} -\frac{1}{t}
\left(\mo{\eps +a}^2-|a|^2\right) (\eps +a)&=&0,\\
\eps\left(1/t_0,x\right)&=&\epsilon_0(x).
\end{array}\right.\end{equation}

\subsection{A priori estimates}

By multiplying the equation (\ref{ee}) by $\overline{\DP{t}\eps}$, and then by taking its real part, we get 
$$\DP{t}\,\frac 12\int|\eps_x (t,x)|^2dx+\frac{1}{t}\int (|\eps(t,x)+a|^2-|a|^2)\Re(\eps(t,x)+a)\overline{\DP{t}\eps}(t,x)dx=0.$$
Therefore we obtain a nice conservation law, that is, 
$$\DP{t}E(t)+\frac{1}{4t^2}\int(|\eps(t,x)+a|^2-|a|^2)^2dx=0,$$
where
$$E(t)=\frac{1}{2}\int|\eps_x(t,x)|^2dx+\frac{1}{4t}\int(|\eps(t,x)+a|^2-|a|^2)^2dx.$$
By integrating from $t'$ to $t$ the energy law, we obtain
$$E(t)+\int_{t'}^t\int(|\eps(\tau,x)+a|^2-|a|^2)^2dx\frac{d\tau}{4\tau^2}=E(t').$$
It follows in particular that for all $t'\leq t$
\begin{equation}\label{grad}
\int|\eps_x(t,x)|^2dx\leq 2E(t'),
\end{equation}
and that
\begin{equation}\label{pot}
\int(|\eps(t,x)+a|^2-|a|^2)^2dx\leq 4t E(t').
\end{equation}

Finally, we shall get a control in time of the mass of $\epsilon$.  By multiplying the equation (\ref{ee}) by $\overline{\eps}$ and taking the imaginary part we obtain the mass law 
\begin{equation}\label{masst}
\DP{t}\,\frac 12\int|\eps (t,x)|^2dx=-\frac{1}{t}\int (|\eps(t,x)+a|^2-|a|^2)\Im a\overline{\eps}(t,x)\,dx.
\end{equation}
By performing a Cauchy-Schwarz inequality in space, we obtain
$$\DP{t}\|\eps (t)\|_2^2\leq \frac{2|a|}{t}\,\|\eps(t)\|_2\,\||\eps(t)+a|^2-|a|^2\|_2.$$
Now we use the upper-bound (\ref{pot}) and get for all $t'\leq t$
$$\DP{t}\|\eps (t)\|_2^2\leq\frac{4|a|\sqrt{E(t')}}{\sqrt{t}}\|\eps(t)\|_2.$$
By integrating from $t'$ to $t$ it follows that for all $t'\leq t$,
\begin{equation}\label{massest}
\|\eps (t)\|_2\leq \|\eps (t')\|_2+c(a,\eps(t'))(\sqrt{t}-\sqrt{t'})\leq c(a,t',\eps(t'))\sqrt{t},
\end{equation}
where the constant depends on $a$, on $t'$, on the energy and $L^2$ norms of $\eps(t')$.

\subsection{Global existence in $H^1$}
First, we solve the equation (\ref{ee}) locally in time. We perform the classical argument of fixed point, done for the cubic equation with constant coefficients (see for example \cite{Gi}, first part of Prop. 4.2.). In the proof, every time that the time $t^{-1}$ appears in the nonlinear terms, we upper-bound it by $t_0$. Therefore we can construct in $X$, the space introduced in \S 2, a solution of (\ref{ee}) living on $(1/t_0,T)$, where $T$ has to verify
$$c(a)\left(T+T^{\frac{1~}{2}}(2\|\epsilon_0\|_2)^2\right)\leq\frac{1}{2}.$$
So it is sufficient that $T$ verifies 
$$T\leq \frac{c(a)}{1+\|\epsilon_0\|_2^4}.$$
In the constant coefficient case, this is enough to infer global existence in $L^2$, since the mass is conserved. In our case, we shall use the control (\ref{massest}).\\ 
Suppose the maximum time of existence is a certain finite $T$. We shall prove that the solution can be defined also after $T$, and therefore the global existence is implied. 
Let $\delta$ a small positive number, to be choosen later. By treating the problem with initial value at time $T-\delta$, we can extend the solution $\eps$ for all time $h$ satisfying
$$h\leq\frac{c(a)}{1+\|\eps(T-\delta)\|_2^4}.$$
By using the control (\ref{massest}), it follows that we can choose such a $h$ verifying furthermore
$$h\geq \frac{c(a,t_0,\epsilon_0)}{1+(T-\delta)^2\|\epsilon_0\|_2^4}.$$
Now, we can choose $\delta$ small enough, such that
$$\frac{c(a,t_0,\epsilon_0)}{1+(T-\delta)^2\|\epsilon_0\|_2^4}>\frac{c(a,t_0,\epsilon_0)}{1+T^2\|\epsilon_0\|_2^4}>\delta,$$
so it follows that $h>\delta$ and so we have extended $\eps$ after the time $T$.\\ 
Therefore we have global existence in $L^2$, provided that the initial value $\epsilon_0$ is finite in energy and $L^2$ norms. By using the Gagliardo-Nirenberg inequalities in the energy definition, it is then enough to have $\epsilon_0$ in ${H}^1$. Finally, in view of the observation (\ref{grad}), the gradient of the solution remains bounded in time, so the global existence is valid in ${H}^1$.\\

\section{A remark on the global existence for the Gross-Pitaevskii equation}

By denoting in (\ref{GP}) $u=\psi-1$, we have
$$\left\{\begin{array}{rcl}
i\partial_t u+\Delta u-\left(|u+1|^2-1\right)(u+1)&=&0,\\
u(0,x)&=&\psi_0(x)-1.
\end{array}\right.$$
We multiply the equation with $\partial_t\overline{u}$, integrate in space and take the real part. We get that the energy
$$E(t)=\frac {1}{2}\int|\nabla u(t,x)|^2dx+\frac {1}{4}\int\left(|u(t,x)+1|^2-1\right)^2dx,$$
is conserved in time. Next, we multiply the equation by $\overline{u}$, integrate in space and take the imaginary part. It follows that
$$\frac{1}{2}\DP{t}\|u(t)\|_2^2=\Im\int\left(|u(t,x)+1|^2-1\right)(u(t,x)+1)\overline{u}(t,x)dx=\Im\int\left(|u(t,x)+1|^2-1\right)\overline{u}(t,x)dx.$$
By performing a Cauchy-Schwarz inequality in space we get 
$$\DP{t}\|u(t)\|_2^2\leq 2\left\||u(t)+1|^2-1\right\|_2\|u(t)\|_2,$$
and so
$$\DP{t}\|u(t)\|_2\leq \left\||u(t)+1|^2-1\right\|_2.$$
The energy conservation allows us to conclude that
$$\DP{t}\|u(t)\|_2\leq 2\sqrt{E(0)}.$$
Therefore we have obtained the claimed a-priori bound on the mass
$$\|u(t)\|_2\leq 2\sqrt{E(0)}\,\,t+\|u_0\|_2,$$
which allows the passage to global existence.

\section{A technical lemma}
\begin{lemma}\label{tech}
Let $x$, $y$ be complex numbers, and $r\geq 0$. Then 
$$\left||x+y|^r-|y|^r\right|\leq c(|y|^{r-1}|x|+|x|^r).$$
Moreover, if $0\leq r\leq1$ then
$$\left||x+y|^r-|y|^r\right|\leq c|y|^{r-1}|x|.$$
Finally, if $x$ is small with respect to $y$, then the last estimate is true for all $r\geq 0$.
\end{lemma}
\begin{proof}
By changing $x=zy$, we have to show for $r\geq 0$, 
$$\left||z+1|^r-1\right|\leq c(|z|+|z|^r),$$
and for $0\leq r\leq 1$, or for $z$ small with respect to $1$,
$$\left||z+1|^r-1\right|\leq c|z|.$$

If $|z|\geq\frac{1}{2}$ then 
$$|z+1|\leq|z|+1\leq 3|z|,$$
and since $r\geq 0$,
$$|z+1|^r\leq 3^r|z|^r.$$
In conclusion,
$$\left||z+1|^r-1\right|\leq 3^r|z|^r\leq c(|z|+|z|^r).$$

If $|z|\leq\frac{1}{2}$ then $0$ is not in the interval $I=[\min\{1,|z+1|\},\max\{1,|z+1|\}]$. We consider the function $f(x)=x^r$ defined on $I$ and we get by the mean value theorem 
\begin{equation}\label{techn}
\left||z+1|^r-1\right|\leq\left||z+1|-1\right|\,r\sup_{\al\in\{1,|z+1|\}}\al^{r-1}\leq |z|r(1+|z+1|^{r-1}).
\end{equation}
If $r\geq 1$, then
$$|z+1|^{r-1}\leq(|z|+1)^{r-1}\leq 2^{r-1},$$
and we get from (\ref{techn}) that
$$\left||z+1|^r-1\right|\leq c|z|.$$
If $r<1$, then 
$$|z+1|\geq 1-|z|\geq\frac{1}{2},$$
and
$$|z+1|^{r-1}\leq 2^{-(r-1)},$$
and it follows again from (\ref{techn}) that
$$\left||z+1|^r-1\right|\leq c|z|.$$
So, under smallness conditions on $|z|$ with respect to $1$, for $r\geq 0$, we have
$$\left||z+1|^r-1\right|\leq c|z|.$$
Otherwise, for general $z$, we get only
$$\left||z+1|^r-1\right|\leq c(|z|+|z|^r).$$\\

Finally, let us treat the case $0\leq r\leq 1$. If $|z+1|\geq\frac{1}{2}$ then
$$|z+1|^{r-1}\leq2^{-(r-1)},$$
and by (\ref{techn}) we get 
$$\left||z+1|^r-1\right|\leq c|z|.$$
If $|z+1|<\frac{1}{2}$ then on one hand $|z+1|^r<2^{-r}<1$, and so
$$\left||z+1|^r-1\right|<1.$$
On the other hand, 
$$|z|\geq 1-|z+1|>\frac{1}{2},$$
and we get again
$$\left||z+1|^r-1\right|<c|z|.$$

\end{proof}

\end{document}